\documentclass[11pt,reqno]{amsart}
\newtheorem{thm}{Theorem}

\newtheorem{prop}{Proposition}

\usepackage{amsfonts}
\usepackage{amsmath}
\usepackage{calligra}
\usepackage{amsmath}
\usepackage{amssymb}

\newcommand{\g}{\mathfrak{g}}
\newcommand{\h}{\mathfrak{h}}
\renewcommand{\k}{{\mathfrak k}}
\renewcommand{\l}{\mathfrak l}

\newcommand{\z}{\mathfrak z}
\renewcommand{\t}{{\mathfrak t}}
\renewcommand{\u}{{\mathfrak u}}

\newcommand{\C}{{\mathbb C}}

\numberwithin{equation}{section}

\begin{document}
\title{Branching laws, some results and new examples }
\author{ Oscar Marquez, Sebastian Simondi,    Jorge A.  Vargas}
\thanks{Partially supported by  CONICET, SECYTUNC
 (Argentina). }
\date{\today }
\keywords{Square integrable representation, admissible restriction, multiplicity formulae }
\subjclass[2010]{Primary 22E46; Secondary 17B10}
\address{ Universidad Federal de Santa María, Santa Maria, RS
Brasil; Facultad de Ciencias Exactas y Natutales, Universidad Nacional de Cuyo, 5500 Mendoza, Argentine; FAMAF-CIEM, Ciudad Universitaria, 5000 C\'ordoba, Argentine}
\email{omskar@gmail.com; ssimondi@uncu.edu.ar; vargas@famaf.unc.edu.ar}

\begin{abstract} For a connected, noncompact matrix simple Lie group $G$ so that a maximal compact subgroup $K$ has three dimensional simple ideal, in this note we analyze the admissibility of the restriction of irreducible square integrable representations for the ambient group when they are restricted to certain subgroups that contains the three dimensional ideal. In this setting we provide a formula for the multiplicity of the irreducible factors. Also, for general $G$ such that $G/K$ is an Hermitian $G$-manifold we give a necessary and sufficient condition so that a square integrable representations of the ambient group is admissible over the semisimple factor of $K.$
\end{abstract}
\maketitle
\markboth{Branching laws}{M-S-V}

\section{Introduction}

Let $G$ be a connected   noncompact simple matrix Lie group. Henceforth, we fix a maximal compact
subgroup $K$ of $G$ and we assume both groups have the same rank. We
also fix $T \subset K$ a maximal torus. Thus, $T$ is a compact Cartan subgroup of
$G$. Under these hypotheses, Harish-Chandra showed there exists irreducible
unitary representations of $G$ so that its matrix coefficients are square integrable
with respect to a Haar measure on $G$. One aim of this note is to write
down explicit branching laws for the restriction of some irreducible square
integrable representation to specific subgroups $H$ of $G$.  A second objective is to show that when  $G$ is simple,  the symmetric space $G/K$ has $G$-invariant quaternionic structure,  and  $H $ is a specific  subgroup   locally isomorphic
to the group $ SU(2, 1) $,  then an irreducible square integrable representations  for $G$  has an admissible restriction to $H$ if and only if it  is a quaternionic discrete series representation. The last objective is to present results on admissible restriction of square integrable representations  to specific subgroups of $G.$  To
begin with, we recall a description of the irreducible square integrable representations
for $G$. Harish-Chandra showed that the set of equivalence classes
of irreducible square integrable representations is parameterized by a lattice
contained in the dual of the Lie algebra of a compact Cartan subgroup.   In order to state our results we need to explicit the parametrization
and set up some notation. As usual, the Lie algebra of a Lie group is denoted
by the corresponding lower case German letter. The complexification
 of a real vector space $V$ is denoted by adding  the subindex $\mathbb C$. However, the root space for a root is denoted
by the real Lie algebra followed by a subindex equal to the root. $V^\star$
denotes the dual space to a vector space $V.$ Let $ \theta $ be the Cartan involution
which corresponds to the subgroup $K,$ the associated Cartan decomposition
is denoted by $\mathfrak g =\mathfrak k +\mathfrak p.$  Let $ \Phi (\mathfrak g, \mathfrak t) $ denote the root system attached to the
Cartan subalgebra $\mathfrak t_\mathbb C.$ Hence, $ \Phi (\mathfrak g, \mathfrak t)  = \Phi_c\cup \Phi_n = \Phi(\k, \t)\cup \Phi_n(\g, \t) $ splits up
as the disjoint  union  of the set of compact roots and the set of noncompact roots. From
now on, we fix a system of positive roots $\Delta$ for $\Phi_c.$ For this note, either the
highest weight or the infinitesimal character of an irreducible representation
of $ K$ is dominant with respect to $\Delta.$ The Killing form on the Lie algebra
$\mathfrak g$ gives rise to an inner product $( \,\, , \,\, )$ in $i\mathfrak t^\star.$ As usual, let $\rho = \rho_g$ denote
half of the sum of the roots for some system of positive roots for $\Phi(\g, \t).$

 A  Harish-Chandra parameter for $G$ is $\lambda \in i\t^\star$ such that $(\lambda, \alpha) \not= 0$
 for every $\alpha \in \Phi (\mathfrak g, \mathfrak t) $, and so that $\lambda+ \rho $ is the differential of a  character of $T.$ To each Harish-Chandra
parameter, $\lambda,$ Harish-Chandra associated a unique irreducible square integrable
representation $(\pi_\lambda^G , V_\lambda)$ of $G$. Moreover, he showed the map $\lambda \mapsto \pi_\lambda^G$ is a bijection from the set of Harish-Chandra parameters dominant with respect to $\Delta$ onto the set of equivalence classes of irreducible square integrable
representations for $G$. For a proof \cite{W1}.

\noindent
In \cite{GW}, the authors have considered quaternionic real form $G$ of a complex
simple Lie group and constructed a specific subgroup $H$ locally isomorphic
to $ SU(2, 1),$ their setting is as follows: a system of positive roots $\Psi$ so that
$ \Delta \subset \Psi \subset  \Phi (\mathfrak g, \mathfrak t) $ is called {\it small} if the maximal root $\beta $ for $\Psi$ is compact,
$\Psi$ has at most two noncompact simple roots $\alpha_0, \alpha_1 $ and after we write $
\beta = n_0\alpha_0+n_1\alpha_1+  $a linear combination of compact simple roots, we have the
inequality $n_0+n_1 \leq 2. $ A noncompact connected simple Lie group $G$ is a {\it quaternionic real form} if $\g$ is an inner form of a complex simple Lie algebra and if a compactly imbedded  Cartan subalgebra $\t$ has the property that  $\Phi(\g,\t)$  admits a small system of positive roots so that $ n_0 + n_1 = 2.$ In \cite{GW},  the  list of the Lie algebras for the quaternionic real forms is presented, we reproduce the list in   Section 3. It can be shown that the set of equivalence classes of the set of quaternionic real forms is equal to the set of equivalence classes of the set of  noncompact simple Lie groups $G$ so that the associated global symmetric space admits a $G$-invariant quaternionic structure.

In order to state Theorem 1, we fix a quaternionic real form $G,$ a small system of positive roots  $\Psi \supset \Delta$
and a noncompact simple root $\alpha$ for $\Psi.$ An irreducible square integrable representation $(\pi_\lambda^G, V)$ is called {\it quaternionic discrete series representation} if the Harish-Chandra parameter $\lambda$ is dominant with respect to $\Psi.$

  For the
quaternionic real form $G$, a particular copy $\h$ of $ \mathfrak{ su}(2, 1)$ contained in  $\g$  is
constructed in  \cite{GW}. For this, they verify the equality
$$ \frac{2(\beta, \alpha)} {
(\alpha, \alpha)}
= 1. $$
Thus, the Lie subalgebra $\h_\mathbb C$ of $\g_\mathbb C$ spanned by the root vectors corresponding
to the roots $\{\pm \alpha,\pm\beta \}$ is isomorphic to $\mathfrak{sl}(3, \C)  $ and invariant under the
conjugation of $\g_\mathbb C$ with respect to $\g.$ Hence, $ \h :=\g \cap \h_\mathbb C$ is a real form for
$\h_\mathbb C$. This real form has a compactly embedded Cartan subalgebra, namely,
$ \u := \t \cap \h.$ Thus, $\h $ is isomorphic to $\mathfrak{su}(2, 1).$ Henceforth, we identify the
set $ \Phi(\h, \u)$ with the subset $\{ \pm \alpha,\pm \beta, \pm (\beta -\alpha)\}$ of $ \Phi (\mathfrak g, \mathfrak t) $. \\ (1.0) Let $H$ denote the
analytic subgroup of $G$ with Lie algebra $\h$. Then,  $L := K \cap H$ is a
maximal compact subgroup for $ H.$ The system $\Phi(\h, \u)$ has three systems of
positive roots to which the root $\beta$ belongs to. The one of our interest is the
non-holomorphic system
$$ \Psi_q := \Psi \cap \Phi(\h, \u) = \{\beta - \alpha, \alpha, \beta \}.$$
The  simple roots for $\Psi_q$
 are $ \beta -\alpha, \alpha. $ For a root $\gamma \in  \Phi (\mathfrak g, \mathfrak t) $, we denote its coroot by
$\check \gamma \in i\t.$   Let $ \Lambda_1,\Lambda_2$ denote the fundamental weights for $\Psi_q,$ labeled so that
$\Lambda_1(\check \alpha) = 0.$

(1.1) Owing to results in \cite{DV}, \cite{W2}, \cite{Kb2}, which we will review in section 2, it
follows that for a Harish-Chandra parameter $\lambda$ dominant with respect to the small system $\Psi$
the irreducible representation $(\pi_\lambda^G
  , V_\lambda)$  restricted to $H$ is an admissible representation.
That is, there exists a sequence of Harish-Chandra parameters for $H$, dominant with respect to $\beta$,  $$ \mu_1, \mu_2, . . . , \mu_j,...\text{ in}\,\, i\u^\star$$
    and there exists  positive integers
  $$ n^{G,H}(\lambda, \mu_j), j = 1, 2, . . .$$
so that the restriction of $ (\pi_\lambda^G
  , V_\lambda)$ to $H$ is unitarily equivalent to the discrete Hilbert sum
$$ \sum_{j=1}^\infty
n^{G,H}(\lambda, \mu_j) \, (\pi_{\mu_j}^H , V_{\mu_j} ).$$ In \cite{GW} it is shown $\Psi_n:=\Psi \cap \Phi_n$ has $2d$ elements.  Our hypothesis that $G$ is a quaternionic real form, forces the
root spaces for the roots  $\pm \beta$ span a three dimensional simple ideal $\mathfrak{su}_2(\beta)$ in $\k$. We
denote by $\k_2$ the complementary ideal to $\mathfrak{su}_2(\beta)$ in $\k.$ Hence, we have the
 decompositions
$$ \t = \mathbb R i \check{ \beta} + (\t \cap \k_2) \,\, \text{and} \,\, \Delta = \{\beta \} \cup \Phi(\k_2, \t \cap \k_2)  \cap \Psi. $$
For each $\lambda \in \mathfrak t_\mathbb C,$ we write $\lambda = \lambda_1 + \lambda_2$ with $\lambda_1 \in \mathbb C \check\beta, \lambda_2 \in \t_{2_\mathbb C} := \mathfrak t_\mathbb C \cap \k_{2_\mathbb C}.$ Let $q_\u : \t^\star \rightarrow \u^\star $ denote the restriction map. \\
(1.3) We will verify (in 2.9) that for   a Harish-Chandra parameter $\lambda$ dominant for the small system $\Psi$, we have $\lambda_2$ is a Harish-Chandra parameter
for $K_2$ or perhaps for a two fold cover of $K_2$. From now on, $\pi_{\lambda_2}^{K_2}$ denotes the irreducible representation for $\k_2$ of infinitesimal character $\lambda_2.$  As usual, $\Delta_{T \cap K_2}\left(\pi_{\lambda_2}^{K_2}\right)$
denotes the set of $T \cap K_2-$weights for the representation $ \pi_{\lambda_2}^{K_2}$
and $ M(\lambda_2, \nu)$
stands for the multiplicity of the weight $\nu \in \Delta_{T\cap K_2} (\pi_{\lambda_2}^{K_2}).$\\
In (2.9) we verify that for $\lambda$ dominant with respect to small system $\Psi$ the weight $\lambda_1+\nu+a\Lambda_1+b q_\u(\lambda_2)$ is dominant with respect to the system $ \Psi_q$ for every $ a, b \in \mathbb Z_{\geq 0},$  and for every $U \cap K_2$-weight $\nu$ of $\pi^{K_2}_{\lambda_2}.$
 One result of this note, is:
 \begin{thm} Let $G$ be a quaternionic real form,  $H$ as in (1.0) and  $(\pi_\lambda^G, V_\lambda)$ a quaternionic discrete series representation for $G$. Then,  $n^{G,H}(\lambda, \mu) \not= 0$ if and only if $\mu= (n+d)\Lambda_1+(m+d)\Lambda_2 +q_\u(\lambda_2) + q_\u(\nu) $ with arbitrary $m,n \in \mathbb{Z}_{\geq 0} $ and $T\cap K_2$-weight $\nu$ for $\pi_{\lambda_2}^{K_2}.$ Moreover

\begin{multline*}
n^{G,H}(\lambda, \mu) = \\
\sum_{\substack{\nu \in \Delta_{T\cap K_2}  \left(\pi_{\lambda_2}^{K_2}
\right),\\
m,n \in \mathbb Z_{\geq 0},\\
\mu=(n+d)\Lambda_1+(m+d)\Lambda_2
+q_\u(\lambda_2)+q_\u(\nu)}}
M(\lambda_2, \nu)
\binom{
m + d - 2}{d - 2}
\binom{n + d - 2}{d - 2}.
\end{multline*}

\end{thm}

\smallskip

A question that naturally arises is: What are the Harish-Chandra parameters for $G,$ dominant with respect to $\Delta, $  so that $\pi_\lambda^G$ has an admissible restriction to $H$? The answer to this question is given in Proposition 1.

\smallskip

A group $G$ locally isomorphic to either $SO(3,n)$   shares with the quaternionic real forms that a  suitable copy of the algebra $\mathfrak{su}_2$ is an ideal in a maximal compactly embedded subalgebra for $\g.$ A group locally isomorphic to $SO(3,2p+1)$ has no square integrable representations. For a group locally isomorphic to $SO(3,2n)$ and $n\geq 2$,  in Proposition 2 we show that no irreducible square integrable representation of $G$   has an admissible restriction to  the usual copy of $"SO(3)"$ contained in $G.$ For the quaternionic group $Sp(1,p)$  the usual factor $"Sp(1)"$ of a maximal compact subgroup is contained in certain  image $H_0$ of $Sp(1,1)$. In  Proposition 3, for a quaternionic discrete series for $Sp(1,p)$,  we show it has admissible restriction to $H_0$,  we compute the Harish-Chandra parameters of the irreducible $H_0$-factors and  their respective multiplicities.

\smallskip

The group $SU(2,1)$ can be mapped into a simple Lie group $G$ in  perhaps several ways by maps $\phi : SU(2,1) \rightarrow G,$ a question is: What are the triple $(G, \pi_\lambda^G, \phi)$ such that $\pi_\lambda$ restricted to the image of $\phi$ is an admissible representation. In \cite{Va1} we find that for the analytic subgroup $H_1$ that corresponds to the  image of $\mathfrak{su}(2,1)$ in the rank one real form of a complex group type $F_4$ no square integrable representation of the ambient group has an admissible restriction to $H_1.$

 We would like to comment that this note grew up from  results in the respective Ph. D. thesis of Sebastian Simondi   and Oscar Marquez  successfully defended at the Faculty of Mathematics Astronomy and Physics at the Universidad Nacional de C\'ordoba, Argentine,  in 2007 and 2011 respectively.

\section{ Proof of Theorem 1}
As in the hypothesis $G$ is a connected, quaternionic simple Lie group and
$H$ is the subgroup locally isomorphic to $SU(2, 1).$
To begin with, we sketch a proof for the statement: {\it For $\lambda$ dominant with respect to the small system $\Psi,$ the representation $(\pi_\lambda^G ,
V_\lambda)$
restricted to $H$ is admissible.} In fact,
for a system of positive roots $ \Sigma \subset \Phi (\mathfrak g, \mathfrak t) $ in \cite{DV} is attached an ideal $ \k_1(\Sigma)$
for the Lie algebra $\k.$ The ideal is equal to the real form of the
ideal of $\k_\mathbb C $ spanned by $\{[Y_\gamma, Y_\phi]: \gamma, \phi \in \Sigma \cap  \Phi_n, Y_\gamma \in  g_\gamma \}$ together with a subspace of
the center $\z_\k$ of $\k$. \\
(2.0) For the system $\Psi, $ cf. \cite{GW} Prop. 1.3, Table 2.5, we have
that any root in $\Phi_c \cap  \Psi$  not equal to $\beta$ is a linear combination of compact
simple roots for $\Psi.$ Thus, for two noncompact roots in $ \Psi ,$ its sum,  is a root
only when the sum is equal to $\beta.$ Thus, $\k_1(\Psi)$ is equal to $\mathfrak{su}_2(\beta)$ plus the contribution of the center. Now, from the list of the quaternionic real forms,
we read that $ \z_\k$ is nonzero only for $G$ locally isomorphic to $SU(2, p).$ For $\mathfrak{su}(2,p)$, in \cite{DV}, it is shown that for $\Psi$ the contribution of $ \z_\k$ to $\k_1(\Psi)$ is just the
zero subspace. Hence, for a quaternionic system $\Psi$ we have
\begin{equation*}
\k_1(\Psi) = \mathfrak{su}_2(\beta).
\end{equation*}

Because of the definition of $H$ we have $ K_1(\Psi)$ is contained in $ H,$ hence, Theorem
1 in \cite{DV} yields that for $\lambda$ dominant with respect to $\Psi$ the representation
$(\pi_\lambda^G, V_\lambda)$ has an admissible restriction to $ H \cap K$ as well as to the subgroup $H.$ In
  \cite{Kb3} we find a different proof of the admissibility.

 \noindent
  Therefore, there exists a sequence of Harish-Chandra parameters
for $ H$, $\mu_1, \ldots , \mu_n, \ldots \in i\u^\star,$ for which, we may assume for every $j,
(\mu_j , \beta) > 0,$ and positive integers $n^{G,H}(\lambda, \mu_j)$ so that the restriction of $(\pi_\lambda^G, V_\lambda)$
restricted to $H $ is equivalent to the Hilbert sum
$$
\sum_j   n^{G,H}(\lambda, \mu_j)\pi_{\mu_j}^H .$$
{\it We are left to compute $\mu_j,$ to show   each $\mu_j$ is dominant for $\Psi_q $ and to compute the
integers} $n^{G,H}(\lambda, \mu_j).$

For this we recall results in \cite{DV}, \cite{He}. For $\gamma \in  i\t^\star.$ (resp in $i\u^\star.$) we consider the
Dirac distribution $\delta_\gamma$ and the the discrete Heaviside distribution defined by
the series
$$
y_\gamma :=
\sum_{n\geq 0}
\delta_{\frac{\gamma}{2}+n\gamma} = \delta_\frac{\gamma}{2}
+ \delta_{\frac{\gamma}{2}+\gamma} + \delta_{\frac{\gamma}{2}+2\gamma} + . . . $$
For any strict multiset $\gamma_1, . . . , \gamma_r$ the convolution $y_{\gamma_1} \star \dots \star y_{\gamma_r}$ is a well
defined distribution. In particular, we have
\begin{equation*} \underbrace{ y_{\gamma} \star \dots \star y_\gamma }_{r} := y_\gamma^r
= \sum_{n \geq 0}
\binom{n + r - 1}{r- 1}
\delta_{(\frac{ r}{2}+n)\gamma}
\end{equation*}
We have
\begin{equation*}
\u_\mathbb C = \mathbb C \check\alpha + \mathbb C \check\beta, \u_\mathbb C \cap \k_{2_\mathbb C} = \u_\mathbb C \cap \t_{2_\mathbb C} = \h_\mathbb C \cap \k_{2_\mathbb C} = \mathbb C( \check\beta - 2 \check\alpha) = \mathbb C \Lambda_2.
\end{equation*}

$ q_\u : \t_\C^\star \rightarrow  \u_\C^\star $ denotes restriction map.

Next, we recall the sub-root system $\Phi_\z := \{\gamma \in \Phi(\k, \t) : q_\u(\gamma) = 0 \}.$
Because of (2.1) and (2.3) $$ \displaystyle{\Phi_\z = \{\gamma : (\gamma, \alpha) = (\gamma, \beta) = 0 \} = \{\gamma \in \Phi(\k_2, \t_2) :
q_{\u\cap \t_2}(\gamma)= 0 \} }. $$ The Weyl group for the system $\Phi_\z$ is denoted by $ W_\z.$ Because
of (1.2) the Weyl group $ W$ for the pair $(\k, \t)$ is equal to the product $ \langle S_\beta \rangle \times
W(\k_2, \t).$ Thus,
$$W_\z \backslash   W =\langle S_\beta \rangle \times W_\z \backslash W(\k_2, \t). $$

Let
 \begin{multline*}
 \Delta(\k/\l) := q_\u[\Psi \cap  \Phi(\k, \t) \backslash \Phi_\mathfrak z] \backslash \Phi(\l, \u) = q_\u[\{\beta \} \cup \Psi(\k_2, \t) \backslash \Phi_\z] \backslash \{\beta\} \\ =
q_{\u \cap \k_2}(\Psi \cap  \Phi(\k_2, t \cap  \k_2) \backslash \Phi_\z) =: \Delta(\k_2/ \u \cap  \k_2).
 \end{multline*}
 We set $\displaystyle{\rho_\z = \frac{ 1}{2} \sum_{ \gamma\in \Psi\cap \Phi_\z } \gamma }$ and
for $\sigma \in i\t^\star$, the Weyl polynomial is defined to be
 $\displaystyle{\varpi (\sigma) :=\frac{\prod_{\gamma \in \Psi	\cap \Phi_\z}
(\sigma,\gamma)}{ \prod_{\gamma \in \Psi	\cap \Phi_\z} (\rho_\z,\gamma) }}. $
As before, we write $\lambda = \lambda_1 + \lambda_2$  with $\displaystyle{\lambda_1 \in \mathbb R i \check\beta}$ and $ \lambda_2 \in \t_2.$  Then, owing to (1.2) for
$\displaystyle{\gamma \in  \Phi(\k_2,\t) \cap  \Psi }$ we have the equality $ \lambda( \check\gamma) = \lambda_2(\check\gamma).$ Thus, $\lambda_2$ is a Harish-
Chandra parameter for $K_2$ or perhaps for a two-fold cover of $K_2$. Actually,
it readily follows that $\lambda_2$ is a Harish-Chandra parameter for $K_2$ if and only if
$\frac{\beta}{2}$ lifts to a character of $T.$   Therefore, if necessary replacing $G$ by a two-fold
cover, we have $\lambda_2$ is a Harish-Chandra parameter for $K_2$.

We now state according to \cite{He} the branching law for the restriction of
the irreducible representation $\displaystyle{\pi_{\lambda_2}^{K_2}}$
of infinitesimal character $\lambda_2$ to the one
dimensional torus $ H \cap K_2 = U \cap K_2.$ The restriction of $\displaystyle{\pi_{\gamma}^{K_2}}$  to $H \cap  K_2$
is the sum of one-dimensional representations $  \sigma_1, . . . , \sigma_r$ with multiplicity
$M(\lambda_2, \sigma_j)$ for $j = 1, . . . , r.$ The formula of Heckman for this particular case
reads
\begin{multline*}
\sum_{ \mu \in \Delta_{U\cap K_2} \left(\pi_{\lambda_2}^{K_2}\right)} M(\lambda_2, \mu) \delta_\mu  \\ =
\sum_{s\in W_\z \backslash W(\k_2, \t)} \epsilon(s) \, \varpi(s\lambda_2) \, \delta_{q_{\u \cap \k_2} (s\lambda_2)} \, \star \, y_{\Delta(\k_2/ \u \cap \k_2)}.
\end{multline*}
Another fact necessary for the proof is a formula in \cite{DV} for the restriction
of $ \displaystyle{\pi_\lambda^G} $ to the subgroup $H.$ The hypothesis for the truth of the formula is
$K_1(\Psi) $ being a subgroup of $H$ which in our case holds because of our choice
of $\Psi$  and $ H.$ The hypothesis on $G$ and on the system $\Psi$ yields for each $w \in W $ the
multiset
\begin{eqnarray*}
S^H_w := [\Delta(\k/\l) \cup q_\u(w\Psi_n)]\backslash \Phi(\h, \u).
\end{eqnarray*}
is strict. This, also follows from an explicit computation of $ S^H_w$ , which, we
will carry out later on. The formula that encodes the parameters $\mu_j$  and the multiplicities $n^{G,H}(\lambda, \mu_j)$ is:
\begin{equation*}
\sum_{\mu \in i\u^\star :(\mu,\beta)>0}
n^{G,H}(\lambda, \mu)(\delta_\mu -\delta_{S_\beta\mu}) =
\sum_{w\in W_\z\backslash W} \epsilon(w)\varpi (w\lambda)\delta_{q_\u(w\lambda) }\star y_{S^H_w}.
\end{equation*}

To elaborate on (2.4) and on (2.6) we recall a few known results. It is
convenient to think of $(\u \cap  \mathfrak{su}_2(\beta))^\star$  (resp.  $\t_2^\star)$ as the linear functionals on $\t$
so that vanishes on $\t_2$ (resp. on $ \u \cap  \mathfrak{su}_2(\beta)),$ hence, for $\lambda_2 \in \t_2,$ we have the
equality $q_\u(\lambda_2) = q_{\u\cap \k_2}(\lambda_2).$
For $ w \in W(\k_2, \t_2)$ we have the equalities
\begin{align*}
q_\u(w\lambda)&= \lambda_1 + q_{\u\cap \k_2}(w\lambda_2) \\
q_\u(wS_\beta\lambda) &= S_\beta(\lambda_1) + q_{\u\cap \k_2}(w\lambda_2)\\
\varpi (w\lambda) &= \varpi (w\lambda_2) \\
\varpi(wS_\beta\lambda) &= \varpi (w\lambda_2).
\end{align*}

From table 2.5 in \cite{GW} it follows that any root in $\Phi(\k_2,\t)$ is linear combination
of compact simple roots in $ \Psi. $ Thus, lemma 3.3 in \cite{HS} yields
\begin{equation*}
 w\Psi_n = \Psi_n \, \text{for} \,  w \in  W(\k_2, \t).
\end{equation*}
In \cite{GW} Proposition 1.3 it is shown that $ \Psi_n = \left\{\gamma \in \Psi :  \frac{2(\beta,\gamma)}{(\beta,\beta)} = 1 \right\},$ and that the
map $\gamma \mapsto \beta - \gamma$ is an involution in $\Psi_n.$ Thus, the number of elements of $\Psi_n$
is an even number $ 2d$ and we may write $$\Psi_n = \{\gamma_2, . . . , \gamma_d, \beta - \gamma_2, . . . , \beta -\gamma_d, \alpha, \beta- \alpha\}.$$ Hence, we have $S_\beta(\Psi_n) = -\Psi_n. $ Also in \cite{GW} Proposition 2 it is shown  that
\begin{equation*}
q_\u(\gamma_j) = \Lambda_1 \, \text{for} \, j = 2, \ldots, d.
\end{equation*}
The equality  $\Lambda_1 + \Lambda_2= \beta $ yields $ q_\u(\beta - \gamma_j) = \Lambda_2. $ From these and (2.8) we
conclude for $ w \in  W(\k_2, \t)$
$$ q_\u(w\Psi_n) = \{ \underbrace{ \Lambda_1, \dots , \Lambda_1}_{d-1}
,\underbrace{\Lambda_2, \dots, \Lambda_2}_{d-1}, \alpha, \beta -\alpha \}. $$

$$ q_\u(wS_\beta\Psi_n) = S_\beta(q_\u(w\Psi_n)) = \{ \underbrace{\Lambda_1, . . . ,\Lambda_1}_{d-1},   \underbrace{\Lambda_2, \dots, \Lambda_2}_{  d-1} ,  \alpha, \alpha -\beta \}.$$
The previous calculations let us conclude.

For $ w \in W(\k_2, \t_2)$,
\begin{align*}
S^H_w &= \{ \underbrace{\Lambda_1, \dots, \Lambda_1}_{d-1}
, \underbrace{ \Lambda_2, \dots, \Lambda_2}_{   d-1} \} \cup \Delta(\k_2/\u \cap  \k_2).\\
 S^H_{S_\beta w} &= \{ \underbrace{-\Lambda_1, ..., -\Lambda_1}_{ d-1}
, \underbrace{-\Lambda_2,\dots, -\Lambda_2}_{  d-1}  \} \cup \Delta(\k_2/\u \cap  \k_2).
\end{align*}
The right hand side of (2.6), after we apply the previous calculations,
becomes equal to
\begin{multline*}
\sum_{s\in W_\z \backslash W(\k_2)}
\epsilon(s) \varpi(s\lambda_2) \delta_{\lambda_1 }\star \delta_{q_{\u\cap \k_2}} (s\lambda_2) \star y^{d-1}_{\Lambda_1} \star y^{d-1}_{\Lambda_2} \star y_{\Delta(\k_2/\u\cap \k_2)} \\
-\sum_{s\in W_\z \backslash W(\k_2) } \epsilon(s) \varpi(s\lambda_2) \delta_{S_\beta\lambda_1 }\star \delta_{q_{\u \cap \k_2 }} (s\lambda_2) \star y_{S_\beta\Lambda_1}^{d-1} \star y_{S_\beta\Lambda_2}^{d-1} \star  y_{\Delta(\k_2/\u\cap \k_2)}
\end{multline*}

$$ =\sum_{ \sigma\in \Delta_{U\cap K_2} (\pi_{ \lambda_2}^{K_2} )}
M(\lambda_2, \sigma) \delta_\sigma \star [\delta_{\lambda_1} \star y^{d-1}_{\Lambda_1} \star y^{d-1}_{\Lambda_2} \star \delta_{S\beta\lambda_1} \star  y^{d-1}_{S_\beta\Lambda_1} \star  y^{d-1}_{S_\beta\Lambda_2}] $$

 \begin{multline*}
 =\sum_{ \sigma,p,q\in \mathbb Z_{\geq 0}} M(\lambda_2, \sigma) \binom{p + d - 2}{d- 2} \binom{q + d - 2}{d - 2} \delta_{\lambda_1+\sigma} \star \delta_{p\Lambda_1} \star \delta_{q\Lambda_2}\\
 + \sum_{\sigma,p,q\in \mathbb Z_{\geq 0}}  M(\lambda_2, \sigma)\binom{p + d -2}{d- 2} \binom{q + d - 2}{d - 2}
\delta_{S_\beta(\lambda_1+\sigma)} \star \delta_{pS_\beta\Lambda_1} \star \delta_{qS_\beta\Lambda_2}.
\end{multline*}

(2.9) We now show: For every $ p, q \in \mathbb Z_{\geq 0},$  and for every $U \cap K_2-$weight $\sigma$ of $\pi^{K_2}_{\lambda_2}$
the weight $\lambda_1+\sigma+p\Lambda_1+q q_\u(\lambda_2)$ is dominant with respect to the system $ \Psi_q= \{\alpha, \beta-\alpha, \beta \}.$

In fact, because of a Theorem of Kostant, every $T_2-$weight of $\displaystyle{\pi^{K_2}_{\lambda_2}}$ lies in the convex hull of $\{s(\lambda_2), s \in W(\k_2, \t) \}.$ Thus, there exists non negative real numbers $c_t$ so that $\sigma =\sum_{t\in W(\k_2,\t)} c_t q_\u(t\lambda_2)$ and  $\displaystyle{\sum_{t} c_t = 1.}$ The hypothesis $\lambda$ is regular and dominant with respect to $\Psi $ yields, $\displaystyle{\lambda_1( \check{ \beta}) = \lambda( \check{ \beta}) > 0.}$

We write $\displaystyle{(\lambda_1 + \sigma + p\Lambda_1 + q\lambda_2, \alpha) = \lambda( \check \beta)\alpha( \check \beta) +\sum_t (q_\u(t\lambda_2), \alpha) + p (\Lambda_1, \alpha)}$ and
$q_\u(t\lambda_2) = (t\lambda_2, \beta -2\alpha)(\beta -2\alpha).$ Now, since $\alpha \in i\u^\star $, we have,
\begin{equation*}
(q_\u(t\lambda_2), \alpha) =
(t\lambda_2, \alpha) = (\lambda_2, t^{-1}\alpha) = (\lambda, t^{-1}\alpha) > 0
\end{equation*} because $t$ is a product of reflections
about compact simple roots for $\displaystyle{\Psi, \alpha \in  \Psi_n}$ and (2.7).

 For
 \begin{align*}
 (\lambda_1 + \sigma+ & p\Lambda_1+q\lambda_2, \beta-\alpha)\\ &=\sum_t c_t(t\lambda_1, \beta-\alpha)+
\sum_t c_t(t\lambda_2, \beta-\alpha)+q(\lambda_2, \beta-\alpha) \\ &=
\sum_t c_t(q_\u( t\lambda_1 + t\lambda_2), \beta - \alpha) + q(\lambda_2, \beta - \alpha) \\ &=
\sum_t c_t(t\lambda, \beta - \alpha) > 0
\end{align*} because of
$\beta -\alpha \in  \Psi_n$, $t \in W(\k_2, \t)$ and $\lambda$ is regular dominant for $\Psi.$
We have concluded the proof of Theorem 1, because we have shown that
the left hand side of (2.6) is expressed as claims the statement of Theorem 1. This finishes the proof of Theorem 1.

{ \bf Note} Wallach in \cite{W2} considered the case the lowest $K$-type for $\pi_\lambda$ is equal to a representation of $\mathfrak{su}_2(\beta)$ times the trivial representation of $K_2.$

\section{Admissible restrictions to $"SU(2,1)"$ of discrete series for quaternionic real forms}
To begin with we list the the Lie algebra of the Lie  groups where Theorem 1 applies. Up to equivalence, the list of the Lie algebras for quaternionic real forms is: $\mathfrak{su}(2, n)$, $\mathfrak{so}(4, n)$, $\displaystyle{EII =
\mathfrak{e}_{6(2)}}$, $EVI = \mathfrak{e}_{7(-5)}$ ,$EIX = \mathfrak{e}_{8(-24)}$, $FI = \mathfrak{f}_{4(4)}$ and $G = \g_{2(2)}.$

For the corresponding groups, we show that a square integrable irreducible
representation for $G$ has an admissible restriction to $ H="SU(2,1)" $ if and only if the
Harish-Chandra parameter is dominant with respect to the small system  $\Psi.$

\begin{prop}
  Let $G$ be a quaternionic real form, a small system of positive roots $\Psi$, $\mathfrak{su}_2(\beta), \k_2, H $ as
in the previous section. Let $\Sigma$ be a system of positive roots in $ \Phi (\mathfrak g, \mathfrak t) $ so
that $\Delta \subset \Sigma$. Then, a square integrable irreducible representation with Harish-
Chandra parameter dominant with respect to $\Sigma$ has an admissible restriction
to $H$ if and only if $\Sigma = \Psi.$
\end{prop}

Proof: From the list of Vogan's diagram, we notice there exists a subgroup
of $H_1$ of $G$ so that $(G, H_1)$ is a symmetric pair and $ H \subset H_1$ and
$T \subset H_1$. Hence, if $\Sigma \supset \Delta$   is a system of positive roots for $ \Phi (\mathfrak g, \mathfrak t) $ so that
some irreducible square integrable representation $\pi_\mu^G$
  with $\mu $ dominant with
respect to $ \sigma$ has admissible restriction to $H,$ then, \cite{Kb3} Theorem 2.8 implies
$\pi_\mu^G$
  has admissible restriction $ H_1$. Owing to \cite{DV} Prop. 2, we have $\k_1(\Sigma)$ is a
subalgebra of $\h_1,$ except for some $G$ locally isomorphic to $ SO(4, 2n)$  the Lie
algebra $ \k$ is the sum of two simple ideals, hence $\k_1(\Sigma)$ is equal to $\mathfrak{su}_2(\beta).$ A  case by case computation forces $\Sigma = \Psi.$ For a group $G$ locally isomorphic
to $SO(4, 2n), n \geq  2 $ we select two different choices of $H_1$  which forces once
again $\k_1(\Sigma) $ to be equal to a copy of $ \mathfrak{su}_2(\beta)$  and $\Sigma = \Psi.$  \qed

\section{Other groups}
A group $G$ locally isomorphic to either $SO(3, n)$ or $Sp(1, n)$ share with
the quaternionic real forms that a copy $\mathfrak{su}_2$ is an ideal in any maximal compactly
embedded subalgebra of   $\g$.  Next, we analyze admissible restriction of square integrable representations to the subgroup corresponding to the copy of $\mathfrak{su}_2$ mentioned in the previous sentence.

 We recall that from a criterium of Harish-Chandra it follows that a group locally isomorphic to $SO(3, 2n + 1)$ has no irreducible  square integrable representation, whereas, a group locally isomorphic
to $SO(3, 2n)$  does have a non empty discrete series.  For a group locally $G$ isomorphic to $SO(3,p)$ a maximally compactly imbedded subalgebra is isomorphic to the direct sum of the ideals $\mathfrak{so}(3), \mathfrak{so}(p).$ For the next statement we denote the analytic subgroup of $G$ corresponding to $\mathfrak{so}(3)$ by $K_1. $

\begin{prop} For a group $G$ locally isomorphic to $SO(3, 2n)$    no irreducible square integrable
representation has an admissible restriction to   $ K_1$.
\end{prop}

Proof: Because, $n \geq 1$  we have that $K_1$ is contained in a subgroup $H_1$
of $G$ locally isomorphic to $ SO(3, 1).$ Next, we recall Theorem 1.2 in \cite{Kb1} which
gives us: if a unitary representation of $G$ has an admissible restriction to
$K_1$ then it has an admissible restriction to $H_1.$ Hence, if irreducible square
integrable representation of $G$ had admissible restriction to $K_1$ we would
have that $H_1$ has a nonempty discrete series, which is not true since $H_1$ is
locally isomorphic to $ SO(3, 1).$ Another proof follows from \cite{DV} and the fact that $K_1(\Psi)$ never is equal to $K_1.$ \qed

For a group $G$ locally isomorphic to $Sp(1, q)$ we fix as maximal compact subgroup $K$ and a compact Cartan subgroup $T.$ Therefore, there exists  an orthogonal basis
$\{ \epsilon_1,  \delta_1, . . . , \delta_q \}$ for $i \t^\star$ and  a system of positive roots $\Sigma$  so that $$\displaystyle{\Sigma \cap  \Phi_c =\{ 2\epsilon_1,  \delta_i \pm \delta_j , 1 \leq  i < j \leq  q, 2\delta_j , j = 1, . . . , q \}}$$ and $\Sigma \cap  \Phi_n = \{ \epsilon_1 \pm \delta_j, j =
1, . . . , q \}.$ The simple roots are $\epsilon_1 - \delta_1$, $\displaystyle{\delta_j - \delta_{j+1}}$, $j = 1, . . . , q, 2\delta_q. $ The
maximal root is $\beta = 2 \epsilon_1.$ It readily follows that $\displaystyle{\k_1(\Sigma) = \mathfrak{su}_2(2 \epsilon_1)}. $  Let $\h_0$
denote the real form of the Lie subalgebra spanned by the root vectors
corresponding to the roots $$ \Phi(\h_0, \u) := \{\pm 2\epsilon_1,\pm 2\delta_1,\pm (\epsilon_1 \pm \delta_1)\}.$$ Then, $\h_0$
is isomorphic to $\mathfrak{sp}(1, 1).$ As for the quaternionic case, let $ H_0 $ denote the
analytic subgroup of $G$ associated to $\h_0.$ Owing to \cite{DV} Theorem 1, we have that
for $\lambda $ dominant with respect to $\Sigma$ the representation $\pi_\lambda^G$
  restricted to $H_0$ is admissible when the Harish-Chandra. Let $\mu_j , n^{G,H_0}(\lambda, \mu_j)$ be as in (1.1). Let
$\Sigma_q:= \Sigma \cap  \Phi(\h_0, \u).$   Let $ HC_{L\cap K_2}\left(\pi^{K_2}_{\lambda_2}\right)$ denotes the set of Harish-Chandra parameters for the
$L\cap K_2-$irreducible factors of the restriction of $\pi^{K_2}_{\lambda_2}$
to the subgroup $L\cap K_2.$ We have,
\begin{prop}Assume $\lambda$ is dominant with respect to $\Sigma$. Then,  for $j = 1, . . .$ the parameters $ \mu_j := \lambda_1 + \sigma + j \epsilon_1 $ are
dominant with respect to $\Sigma_q.$ Besides,  $  n^{G,H_0}(\lambda, \mu) \not= 0 $ if and only if $\mu =\mu_j$ for some $j.$ Moreover
\begin{equation*}
n^{G,H_0}(\lambda, \mu) =\sum_{\substack{\sigma\in HC_{L\cap K_2} (\pi_{\lambda_2}^{K_2}), \\ p \in \mathbb Z_{\geq 0} \\ \mu=\lambda_1+\sigma+p \epsilon_1}}
M(\lambda_2, \sigma) \binom{p + 2q - 3}{2q - 3}.
\end{equation*}
\end{prop}

Proof: We begin writing the equalities (2.4) and (2.6) for the setting of
the Proposition. For this particular case (2.4) reads
\begin{multline*}
\sum_{ \mu\in HC_{L\cap K_2}( \pi^{K_2}_{\lambda_2})} M(\lambda_2, \mu) \sum_{r\in W(L\cap K_2, U\cap K_2)}
\epsilon(r)\delta_{r\mu} \\ =\sum_{s\in W_\z \backslash W(\k_2,\t)}
\epsilon(s) \varpi(s\lambda_2) \delta_{q_{\u \cap \k_2 }} (s\lambda_2) \star y_{\Delta(\k_2/ \l\cap \k_2)}.
\end{multline*}
The multiset
\begin{equation*}
S^{H_0}_w := [\Delta(\k/\l) \cup q_\u(w\Psi_n)]\backslash \Phi(\h, \u)
\end{equation*}
is strict. This, follows from an explicit computation of $ S^{H_0}_w,$ which, we will
carry out after the next formula. The formula  (2.6)  becomes:
\begin{multline*}
\sum_{\substack{\mu\in i\u^\star \\
(\mu, \epsilon_1)>0, \,(\mu,\delta_1)>0} }n^{G,H}(\lambda, \mu)\left(\sum_{t\in W(L,U)} \epsilon(t) \delta_{t\mu}\right) $$
\\ =
\sum_{w\in W_\z \backslash W} \epsilon_1(w) \varpi(w\lambda) \delta_{q_\u(w\lambda)} \star y_{S^H_w}.
\end{multline*}

In this case $\u\cap \k_1 = \mathbb Ri \check \epsilon_1$, $\u\cap \k_2 = \mathbb R i \check \delta_1$,  $\l \cap \k_2 = \mathfrak{su}_2(2\delta_1)$. Furthermore
\begin{eqnarray*}
\Psi\cap \Phi_\z = \{\delta_i\pm \delta_j , 2 \leq i < j \leq q \}, \, \,
\end{eqnarray*}  $W = \langle S_{2\epsilon_1} \rangle \times W(\k_2, \t_2)$  and $W_\z\backslash W = \langle S_{2\epsilon_1} \rangle \times W_\z \backslash W(\k_2, \t_2).$

For $ w \in W(\k_2, \t_2)$,   $w\Psi_n = \Psi_n$, $wS_{2\epsilon_1} \Psi_n = -\Psi_n.$ and

\begin{align*}
q_\u(w\Psi_n) \backslash \Phi(\h, \u) &= \underbrace{ \{\epsilon_1, . . .,\epsilon_1 \} }_{2(q-1)}.\\
q_\u(wS_{2\epsilon_1} \Psi_n) \backslash \Phi(\h, \u) &= \underbrace{ \{-\epsilon_1, . . . , -\epsilon_2 \} }_{2(q-1)}.\\
\Delta(\k, \l) = q_\u(\Psi_c \backslash \Phi_\z) \backslash \Phi(\h, \u) &=  \{2\epsilon_1, 2\delta_1, \underbrace{ \delta_1, . . . , \delta_1}_{2(q-1)} \} \backslash \Phi(\h, \u) = \Delta(\k_2, \l\cap \k_2).
\end{align*}

Therefore, for $ w \in W(\k_2, \t) $ we have,

\begin{align*}
S^{H_0}_w &=\underbrace { \{\epsilon_1, . . . , \epsilon_1\} }_{2(q-1)} \cup \Delta(\k_2, \l \cap  \k_2).\\
 S_{S_2 \epsilon_1 w}^{H_0}  &= \underbrace{ \{ \epsilon_1, ... , \epsilon_1 \} }_{2(q-1)} \cup \Delta(\k_2, \l \cap  \k_2).
\end{align*}

After replacing $S^{H_0}_w$ by the result obtained in the previous line, the right
hand side of the formula similar to the  one in (2.6) becomes

 \begin{eqnarray*}
 \sum_{\substack{t\in \{1,S_{2\epsilon_1} \} \\ s \in W_\z \backslash W(\k_2,\t_2)}} \epsilon(t)\epsilon(s) \varpi(s\lambda_2) \delta_{t\lambda_1} \star \delta_{q_\u(s\lambda_2)} \star   y_{\Delta(\k_2/ \l\cap  \k_2)} \star y^{2(q-1)}_{t \epsilon_1}\\
 \end{eqnarray*}
\begin{align*}
 &=\sum_{t\in \{1,S_{2\epsilon_1} \} }  \epsilon(t) \sum_{\substack{r\in W(\k_2,\t_2) \\ \sigma\in  HC}} \epsilon(r) M(\lambda_2, \sigma) \delta_{t\lambda_1+r\sigma } \star y^{2(q-1)}_{t \epsilon_1}\\
&=\sum_{\substack{t,r,\sigma \\ p\in \mathbb Z_{\geq 0}}} \epsilon(t) \epsilon (s)M(\lambda_2, \sigma) \binom{p + 2(q-1) - 1}{2(q - 1) - 1}  \delta_{tr(\lambda_1+\sigma+p \epsilon_1)} \\
&=\sum_{p \geq 0} \sum_{\sigma} M(\lambda_2, \sigma) \binom{p + 2q - 3}{2q - 3}  \sum_{w \in W(L,U)} \epsilon(w)\delta_{w(\lambda_1+\sigma+p \epsilon_1)}.
\end{align*}
By a reasoning similar to (2.9) we obtain that $\lambda_1+\sigma+p \epsilon_1$ is dominant with
respect to $\Psi_q $ and we have concluded the proof of  Proposition 3. \qed

\section{  Simondi's Thesis}

The Ph. D. thesis of Sebastian Simondi was defended by March 2007. Most of his results were verified in a case by case checking. By now, some of his results are a consequence of work of Toshi Kobayashi \cite{Kb2}, \cite{Kb3} and \cite{DV}. This section does not include proofs of the stated results, we will point out those results that follows from the work of  T. Kobayashi, Kobayashi-Oshima, DV and those results we  believe are still not in the literature. A complete version of Simondi's thesis is on:\\  www2.famaf.unc.edu.ar/publicaciones/documents/seried/DMat49.pdf.

(5.1) We now describe the results. For this, we fix a noncompact connected matrix simple Lie group $G$ a maximal compact subgroup $K$ for $G.$ Henceforth,  $H$ is a closed  reductive subgroup of $G$ so that $L:=H\cap K$ is a maximal compact subgroup for $H$ and  that $(G,H)$ is a symmetric pair. Hence, $(K,L)$ is Riemannian symmetric pair.

As in the previous sections we assume $G$ admits irreducible square integrable representations, we would like to point out that  in the course of the computation was made an extensive use of the description for the set of equivalence classes of square integrable  irreducible representations given by  Harish-Chandra in terms of Harish-Chandra parameters.

In \cite{DV}, \cite{Kb1}, \cite{Kb2} and \cite{Kb3} we find criteria for checking whether or not the restriction of an square integrable representation for $G$ is an admissible representation for a subgroup $H.$ By mean of these criteria, the classification of the symmetric pairs given by Berger and a case by case checking,  we have,

\begin{thm} Assume $(G,H)$ is a symmetric pair and $(\pi, V)$ is  an irreducible square integrable representation for $G$. If $\k$ is a simple Lie algebra, then the restriction of $\pi$ to $H$ is not an admissible representation.
\end{thm}

Nowadays, this result follows from \cite{DV} or from the work of \cite{KO}.

\smallskip
For the next result we fix a maximal compact connected subgroup $L^\prime$ for $K$ so that the rank of $K$ is equal to the rank of $L^\prime.$

\begin{thm} Let $(\pi, V)$ be an irreducible square integrable representation for $G$. We assume $\k$ is a simple Lie algebra. Then, $\pi$ restricted to $L^\prime$ is not an admissible representation.
\end{thm}

When $L^\prime$ is a maximal compact subgroup of a reductive subgroup $H$ of $G$ so that $(G,H)$ is symmetric pair  Theorem 3 follows from Theorem 2 and results in \cite{DV}.  For the other subgroups $L^\prime$ the proof has been done in a case by case checking based on the classification of the equal rank maximal subgroups of $K$ obtained by Borel-de Siebenthal and work of Toshi Kobayashi on criteria on admissibility of restriction of representations.

\smallskip

\noindent
Under the hypothesis $\k$ is not a simple Lie algebra, $(G,H)$ a symmetric pair and the subgroups $L,K$ are of the same rank, we obtain a complete list, in the language of Harish-Chandra parameters, of the square integrable representations for $G$ which do not have an admissible restriction to $H.$ Nowadays this results are included in \cite{KO} \cite{Va2}.

\smallskip

For the last result of this note  we further  assume $(G,K)$ is an Hermitian symmetric  pair.  Then, the center of $K$ is a one dimensional torus. Let $K_{ss}$ denote the semisimple factor of $K.$ We fix a maximal torus $T$ for $K.$  The hypothesis on $(G,K)$ allows us to choose, once for all, a holomorphic system of positive roots $\Psi_h$ in $\Phi(\g, \t)$. In \cite{Kb1} it is shown that either a holomorphic or a  antiholomorphic discrete series for $G$ has an admissible restriction to $K_{ss}$ if and only if  $G/K$ is not a tube domain. The next result gives a criteria which allows to determine when an  arbitrary irreducible square integrable representation  has  admissible restriction to $K_{ss}$.  For this we recall set of equivalence classes for  irreducible square integrable representations is parameterized by the set of  Harish-Chandra parameters $\lambda$  dominant with respect to $\Psi_h \cap \Phi_c$. The regularity of   $\lambda$  determines a system of positive roots $\Psi_\lambda:=\{ \alpha \in \Phi( \g, \t): \lambda (\check{\alpha})>0 \}$ which satisfies  $\Psi_\lambda \cap \Phi_c = \Psi_h \cap \Phi_c. $ In Table 1 we list for each Hermitian symmetric  pair $(G,K)$ subsets $I,\tilde I$ of $\Psi_h.$
\begin{thm}Assume $(G,K)$ is an Hermitian symmetric pair and fix an  irreducible square integrable representation $(\pi_\lambda^G, V_\lambda)$  for $G$ of Harish-Chandra parameter $\lambda$ dominant with respect $\Psi_h \cap \Phi_c$. Then, $\pi_\lambda^G$  restricted to $K_{ss}$ is an admissible representation if and only if either $I$ or $\tilde I$ is a subset of $\Psi_\lambda.$
\end{thm}

 The proof of the last Theorem is carried out in basis of classification of Hermitian symmetric  pairs and   criteria due to Kobayashi \cite{Kb2}, \cite{Kb3}.
\begin{center}
\begin{tabular}{|c|c|l|}
\hline
 $\frak{g} $ & $ \frak{k} $ & \\
 \hline
 $\frak{sp}(n, \mathbb{R}) $ & $ \frak{u}(n) $ &
 \begin{tabular}{l}
 \medskip
If $n=2l$ \\
 $I=\{(e_k+e_{n-k+1})\}_{k=1}^{l},$ \\
\medskip
 $ \tilde{I}=-I$  \\
 If $n=2l+1$\\
 $I=\{2e_{l+1}\} \cup \{(e_k+e_{n-k+1})\}_{k=1}^{l},$ \\
\medskip
 $ \tilde{I}=-I$
  \\
 \end{tabular}\\
 \hline
 $\frak{so}^{*}(2n) $& $ \frak{u}(n)$ & \begin{tabular}{l}
If $n=2l$\\
 $I=\{(e_k+e_{n-k+1})\}_{k=1}^{l},$ \\
 $ \tilde{I}=-I$  \\
If $n=2l+1$\\
$I=\{(e_k+e_{n-k+1})\}_{k=1}^{l} \cup \{ e_{l+1}+e_{l+2} \},$  \\
 $ \tilde{I}=\{(-e_k-e_{n-k+1})\}_{k=1}^{l} \cup \{ -e_{l}-e_{l+1} \}$ \\
 \end{tabular}\\
 \hline
 $\frak{su}(p,q)$& $\frak{su}(p)\oplus \frak{u}(q)$ &
 \begin{tabular}{l}
 $I=\left\{e_i-e_{\gamma_i} \right\}_{i=1}^{p}$ \\
\\
 $\tilde{I}=\left\{-e_i+e_{b_i} \right\}_{i=1}^{p}$
 \end{tabular} \\
 \hline
$\frak{e}_{6(-14)}$ & $ \frak{so}(10) \oplus \frak{so}(2)$ &
\begin{tabular} {l}
$I=\{ \varepsilon_1 , \varepsilon_2,e_1+e_5, e_2+e_5\}$\\
$\tilde{I} = - I$
\end{tabular}
\\ \hline
$\frak{e}_{7(-25)}$ & $ \frak{e}_6 \oplus \frak{so}(2)$ &
\begin{tabular}{l} $I= \{\eta_1, \eta_2,e_1+e_6 \}$
\\ $\tilde{I} = -I$ \end{tabular}
\\
\hline
\end{tabular}
\end{center}
\begin{center}
Table 1
\end{center}

Here, in  notation of Bourbaki:

\begin{tabular}{l}
\bigskip
$\gamma_i= i+\left[\frac{(i-1)(q-p)}{p}\right]+p,$ \text{for} $1 \leq i \leq p,$\\
\bigskip
$ b_i=\left\{
\begin{array}{l}
  i+1+[\frac{i(q-p)}{p}]+p
  \quad\quad\quad \text{if} \,\, \frac{i(q-p)}{p} \notin \mathbb{Z},  \\
  i+\frac{i(q-p)}{p}+p \,\,
  \quad\quad\quad \quad\quad \text{if} \,\, \frac{i(q-p)}{p} \in \mathbb{Z}, \\
\end{array}%
\right. $ \text{for} $1 \leq i \leq p,$ \\
\bigskip
$\varepsilon_1=\frac{1}{2}(-e_1-e_2-e_3-e_4+e_5-e_6-e_7+e_8)$, \\
\bigskip
$\varepsilon_2=\frac{1}{2}(-e_1-e_2+e_3+e_4+e_5-e_6-e_7+e_8)$, \\
\bigskip
$\eta_1=\frac{1}{2}(-e_1+e_2-e_3-e_4+e_5+e_6-e_7+e_8)$ \\
\bigskip
$\eta_2=\frac{1}{2}(-e_1-e_2+e_3+e_4-e_5+e_6-e_7+e_8)$ \\
\end{tabular}

\providecommand{\MR}{\relax\ifhmode\unskip\space\fi MR }
\providecommand{\MRhref}[2]{%
  \href{http://www.ams.org/mathscinet-getitem?mr=#1}{#2}
}
\providecommand{\href}[2]{#2}

\end{document}